\documentclass[a4paper]{amsart}
\pdfoutput=1
\usepackage{amssymb,color}
\usepackage{graphicx}
\usepackage{upref}
\definecolor{darkred}{rgb}{0.6,0,0}
\definecolor{darkgreen}{rgb}{0,0.5,0}
\definecolor{darkmagenta}{rgb}{0.5,0,0.5}
\usepackage[pdftex, colorlinks=true,
  linkcolor=darkred, citecolor=darkgreen, urlcolor=darkmagenta]{hyperref}

\title[]{Addendum to ``The Kolmogorov--Riesz compactness theorem''  (Expo. Math. 28:385--394 (2010))}

\author[Hanche-Olsen]{Harald Hanche-Olsen}
\address[Hanche-Olsen]{\newline
    Department of Mathematical Sciences,
    Norwegian University of Science and Technology,
    NO--7491 Trondheim, Norway }
\email[]{\href{hanche@math.ntnu.no}{hanche@math.ntnu.no}}
\urladdr{\href{http://www.math.ntnu.no/~hanche}{www.math.ntnu.no/\~{}hanche}}

\author[Holden]{Helge Holden}
\address[Holden]{\newline
    Department of Mathematical Sciences,
    Norwegian University of Science and Technology,
    NO--7491 Trondheim, Norway}
\email[]{\href{holden@math.ntnu.no}{holden@math.ntnu.no}}
\urladdr{\href{http://www.math.ntnu.no/~holden}{www.math.ntnu.no/\~{}holden}}

\subjclass[2010]{Primary: 46E30, 46E35; Secondary: 46N20}

\keywords{Kolmogorov--Riesz' theorem}

\thanks{Supported in part by the Research Council of Norway.}

\newcommand{\RR}{\mathbb R}
\newtheorem*{theorem}{Theorem}
\theoremstyle{remark}
\newtheorem{remark}{Remark}
\newcommand{\withdelims}[4]{\mathopen#1#2#3\mathclose#1#4}
\newcommand{\defdelims}[4][]
  {\newcommand{#2}[2][#1]{\withdelims{##1}{#3}{##2}{#4}}}
\defdelims{\abs}||
\defdelims{\norm}\|\|
\defdelims{\parens}()
\defdelims[\big]{\setof}\{\}
\newcommand{\dott}{\, \cdot\,}

\allowdisplaybreaks

\begin{document}

\begin{abstract} 
We show how to improve on Theorem 10 in \cite{HOH}, describing when subsets in 
$W^{1,p}(\RR^n)$ are totally bounded subsets of $L^q(\RR^n)$ for $p<n$ and $p\le q<p^*$. This improvement was first shown in \cite{DFS} in the context of Morrey--Sobolov spaces. 
\end{abstract}

\maketitle

\noindent
We show the following improvement of Theorem 10 in \cite{HOH}:
\begin{theorem}
Assume $1\le p<n$ and $p\le q<p^*$, where
\[
  \frac1{p^*}=\frac1p-\frac1n,
\]
and let $\mathcal{F}$ be a bounded subset of $W^{1,p}(\RR^n)$.
Assume that for every $\epsilon>0$ there exists some $R$ so that, for every
$f\in\mathcal{F}$,
\begin{equation} \label{eq:bet}
 \int_{\abs{x}>R}\abs{f(x)}^p\,dx<\epsilon^p.
\end{equation}
Then $\mathcal{F}$ is a totally bounded subset of $L^q(\RR^n)$.
\end{theorem}

\begin{remark} The improvement over the original consists of replacing
\begin{equation*}
 \int_{\abs{x}>R}\bigl(\abs{f(x)}^p+\abs{\nabla f(x)}_p^p\bigr)\,dx<\epsilon^p
\end{equation*}
by the weaker inequality \eqref{eq:bet}.
(Here $\abs{\dott}_p$ is the $l^p$ norm on $\RR^n$.)
\end{remark}
\begin{remark} The improvement we prove here was first shown in the more complicated setting of Morrey--Sobolev spaces in \cite{DFS}. Here we present a direct argument in the setting of \cite{HOH}.
\end{remark}

\begin{proof}
  The Sobolev embedding theorem states that
  $W^{1,p}(\RR^n)\subset L^q(\RR^n)$
  and that the inclusion map is bounded.
  Thus $\mathcal{F}$ is a bounded subset of $L^q(\RR^n)$.

  We will need the Gagliardo--Nirenberg--Sobolev inequality, which
  states that there is a constant $C$ (only dependent on $p$ and $n$)
  such that
  \begin{equation*}
    \norm{f}_{p^*}\le C\norm{\nabla f}_p
  \end{equation*}
  for all $f\in C_c^1(\RR^n)$.
  For a proof, see, e.g., \cite[Sect.~5.6.1, p.~263]{evans}.
  This inequality extends to any $f\in W^{1,p}(\RR)$:
  Let $f_n\in C_c^1(\RR)$ converge to $f$ in $W^{1,p}(\RR)$.
  Then the inequality implies that $(f_n)$ is Cauchy in $L^{p^*}(\RR)$,
  so it has a limit, which must be $f$ itself,
  since some subsequence converges pointwise.
  The continuity of the norms finishes the argument.

  The above inequality and the interpolation inequality
  $\norm{f}_q \le\norm{f}_p^\theta\, \norm{f}_{p^*}^{1-\theta}$
  where
  \[
    \frac{1}{q}=\frac{\theta}{p}+\frac{1-\theta}{p^*}\qquad(0<\theta\le1)
  \]
  yield
  \[
    \norm{f}_q
    \le C^{1-\theta}\norm{f}_p^\theta\,\norm{\nabla f}_{p}^{1-\theta}.
  \]

  Now let $\epsilon>0$, and pick $R$ as in the statement of the theorem.
  Let $\phi\in C_c^\infty(\RR^n)$ be a function
  with $0\le\phi\le1$ and $\abs{\nabla\phi}_p\le1$
  satisfying $\phi(x)=1$ when $\abs{x}\le R$.

  Then $\phi\mathcal{F}=\setof{\phi f\colon f\in\mathcal F}$ is
  bounded in $W^{1,p}(\RR^n)$, and by \cite[Theorem 10]{HOH}
  (or the usual Rellich--Kondrachov theorem on a ball of radius $R+2$),
  $\phi\mathcal{F}$ is totally bounded in $L^q(\RR^n)$.

  We find that every $f\in\mathcal{F}$ is uniformly approximated in
  the $L^q$ norm by $\phi f$:
  \begin{align*}
    \norm{f-\phi f}_q
    &\le C^{1-\theta}\norm{(1-\phi)f}_p^\theta\,\norm{\nabla((1-\phi)f)}_{p}^{1-\theta}
  \\&\le C^{1-\theta}\epsilon^\theta\,
      \norm{(1-\phi)\nabla f-f\nabla\phi}_{p}^{1-\theta}
     \\&\le C^{1-\theta}\epsilon^\theta\,
     \parens[\big]{\norm{(1-\phi)\nabla f}_p+\norm{f\nabla\phi}_p}^{1-\theta}
    \\&\le\parens[\big]{C2^{1-1/p}\norm{f}_{1,p}}
          ^{1-\theta}
        \epsilon^\theta
  \\&\le M\epsilon^\theta
  \end{align*}
  where the constant $M$ depends only on $n$, $p$, $q$, and $\mathcal F$.
  In the next to last line, we used Jensen's inequality on the form
    $u+v\le2^{1-1/p}(u^p+v^p)^{1/p}$ when $u,v\ge0$,
    together with the definition of the $W^{1,p}$ norm.

  Thus every member of $\mathcal{F}$ lies within a distance $M\epsilon^\theta$
  of $\phi\mathcal{F}$ in $L^q$ norm.
  Since $\phi\mathcal{F}$ is totally bounded and $M\epsilon^\theta$
  can be made arbitrarily small,
  it follows that $\mathcal{F}$ is totally bounded.
\end{proof}

\end{document}